# ON 2-MOVABLE DOMINATION IN THE JOIN AND CORONA OF GRAPHS

**Ariel C. Pedrano[1] and Rolando N. Paluga[2]**

[1]Department of Mathematics and Statistics
University of Southeastern Philippines
Davao City, Philippines
e-mail: ariel.pedrano@usep.edu.ph

[2]Mathematics Department
Caraga State University
Butuan City, Philippines
e-mail: rnpaluga@carsu.edu.ph

## Abstract

Let $G$ be a connected graph. Then a non-empty $S \subseteq V(G)$ is a 2-movable dominating set of $G$ if $S$ is a dominating set and for every pair $x, y \in S$, $S - \{x, y\}$ is a dominating set in $G$, or there exist $u, v \in V(G) \setminus S$ such that $u$ and $v$ are adjacent to $x$ and $y$, respectively, and $(S \setminus \{x, y\}) \cup \{u, v\}$ is a dominating set in $G$. The 2-movable domination number of $G$, denoted by $\gamma_m^2(G)$, is the minimum cardinality of a 2-movable dominating set of $G$. A 2-movable dominating set with cardinality equal to $\gamma_m^2(G)$ is called $\gamma_m^2$-set of $G$.

This paper obtains 2-movable domination numbers for the corona and join of graphs.

## 1. Introduction

All graphs considered in this paper are all connected, finite, simple and undirected. Let $G = (V, E)$ be a connected, finite, simple and undirected graph. The graph $G$ has a vertex set $V = V(G)$ and an edge set $E = E(G)$. Further, let the order of the graph $G$ be $p$, that is, $|V| = |V(G)| = p$ and the size of the graph $G$ be $q$, that is, $|E| = |E(G)| = q$.

One of the interesting fields of Graph Theory is graph domination. Ore [2] introduced the terms "dominating set" and "domination number" in his book entitled *Theory of Graphs* which was published in 1962. A subset $S$ of $V(G)$ is a *dominating set* of $G$ if for every $v \in V(G) \backslash S$, there exists $u \in S$ such that $uv \in E(G)$, that is, the closed neighbourhood of $S$ is the vertex set of $G$. The domination number of $G$ is denoted by $\gamma(G)$ which refers to the smallest cardinality of a dominating set of $G$. A dominating set of $G$ with cardinality equal to $\gamma(G)$ is called a $\gamma$-*set* of $G$.

Blair et al. [1] introduced a new variant of domination and called it 1-*movable domination*. A 1-movable dominating set is a dominating set $S \subseteq G$ such that for every $v \in S$, at least one of the following two conditions holds, i.e., $S - \{v\}$ is a dominating set, or there exists a vertex $u \in (V(G) - S) \cap N(v)$ such that $(S - \{v\}) \cap \{u\}$ is a dominating set. Inspired by the study of Blair et al. [1], we are interested in exploring the concept of 2-movable domination in the join and corona of graphs.

## 2. Basic Concepts

**Definition 2.1.** Let $G$ be a connected graph. A non-empty $S \subseteq V(G)$ is a 2-*movable dominating set* of $G$ if $S$ is a dominating set and for every pair $x, y \in S$, $S - \{x, y\}$ is a dominating set in $G$, or there exist $u, v \in V(G) \backslash S$ such that $u$ and $v$ are adjacent to $x$ and $y$, respectively, and $(S \backslash \{x, y\}) \cup \{u, v\}$ is a dominating set in $G$. The 2-*movable domination number* of $G$, denoted by $\gamma_m^2(G)$, is the minimum cardinality of a 2-movable dominating set of $G$. A 2-movable dominating set with cardinality equal to $\gamma_m^2(G)$ is called a $\gamma_m^2$-*set* of $G$.

Consider the graph $G$ in Figure 1. The set $S_1 = \{v_2, v_7, v_9\}$ is a $\gamma$-set, $S_2 = \{v_2, v_6, v_7, v_9\}$ is a $\gamma_m^1$-set and $S_3 = \{v_2, v_5, v_6, v_7, v_{12}\}$ is a $\gamma_m^2$-set. Thus, $\gamma(G) = |S_1| = 3$, $\gamma_m^1(G) = |S_2| = 4$ and $\gamma_m^2(G) = |S_3| = 5$.

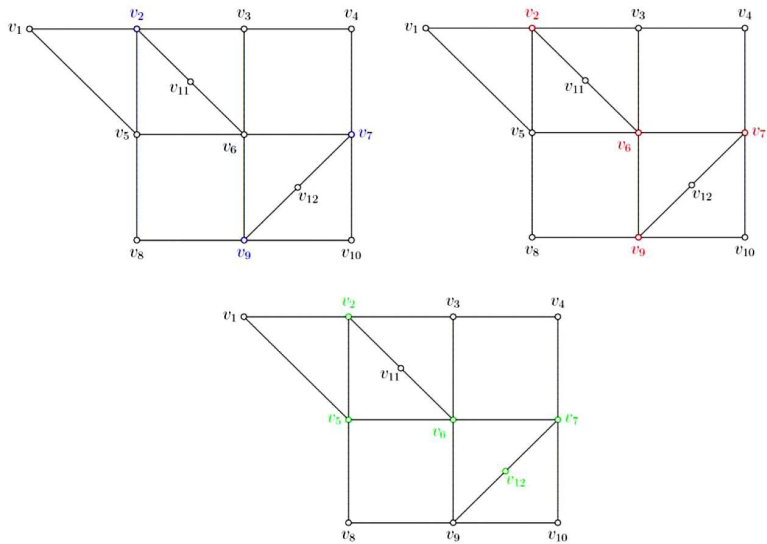

**Figure 1.** The graph $G$.

### 3. Results

The next remark follows directly from Definition 2.1.

**Remark 3.1.** For any connected graph $G$ of order $n \geq 4$,

$$\gamma_m^2(G) \geq 2.$$

**Theorem 3.2.** *For any connected graph $G$ of order $n \geq 4$,*

$$\gamma(G) \leq \gamma_m^1(G) \quad \text{and} \quad \gamma(G) \leq \gamma_m^2(G).$$

**Proof.** Suppose $S$ is a 1-movable dominating set of $G$. It follows that $S$ is a dominating set of $G$. Therefore, $\gamma(G) \leq \gamma_m^1(G)$.    □

The next theorem provides the exact value of 2-movable domination number in the join of two graphs.

**Theorem 3.3.** *Let $G$ and $H$ be graphs of order at least 2. Then*

$$\gamma_m^2(G + H) = 2.$$

**Proof.** Let $|V(G)| \geq 2$ and $|V(H)| \geq 2$, where $V(G + H) = V(G) \cup V(H)$ and $|V(G + H)| \geq 4$. Let $S = \{u, v\}$, where $u \in V(G)$ and $v \in V(H)$. We want to show that $S$ is a dominating set of $G + H$. Observe that $V(H) \subseteq N[u]$ and $V(G) \subseteq N[v]$. This implies that $V(G + H) = V(G) \cup V(H) \subseteq N[u] \cup N[v] = N[S]$. Hence, $S$ is a dominating set of $G + H$.

Now, since $|V(G)| \geq 2$ and $|V(H)| \geq 2$, there exist $u' \in V(G)$ and $v' \in V(H)$ such that $u \neq u'$ and $v \neq v'$. Thus, $(S \backslash \{u, v\}) \cup \{u', v'\} = \{u', v'\}$. By following the same argument above, $\{u', v'\}$ is a dominating set of $G + H$. Moreover, $u' \in N[v]$ and $v' \in N[u]$. Hence, $S$ is a 2-movable dominating set of $G + H$. Therefore, $\gamma_m^2(G + H) \leq 2$. By Remark 3.1, $\gamma_m^2(G + H) = 2$.    □

**Lemma 3.4.** *Let $T$ be a dominating set of $G \circ H$ and $a \in V(G)$. If $a \notin T$, then $S_a = T \cap V(H^a)$ is a dominating set of $H^a$.*

**Lemma 3.5.** *For any corona graph $G \circ H$ and $a \in V(G)$, if $T$ is a 2-movable dominating set of $G \circ H$, $T_a = T \cap V(a + H^a)$ and $u \in S_a$, then*

(i) $T_a \setminus \{a, u\}$ *is a dominating set of* $H^a$, *or*

(ii) *there exist $x_a, x_u \in V(a + H^a) \setminus T_a$ such that $ax_a, ux_u \in E(a + H^a)$ and $(T_a \setminus \{a, u\}) \cup \{x_a, x_u\}$ is a dominating set of $H^a$,*

(iii) *there exists $x_u \in V(a + H^a) \setminus T_a$ such that $ux_u \in E(a + H^a)$ and $(T_a \setminus \{a, u\}) \cup \{x_u\}$ is a dominating set of $H^a$.*

The next theorem provides the exact value of 2-movable domination number in the corona of two graphs.

**Theorem 3.6.** *Let $G$ and $H$ be connected graphs such that $|V(G \circ H)| \geq 4$. Then*

$$\gamma_m^2(G \circ H) = |V(G)|\gamma(H).$$

**Proof.** For each $x \in V(G)$, let $H^x$ be a copy of $H$ corresponding to vertex $x$. Further, suppose that $S_x$ is a minimum dominating set of $H^x$. Let $S = \bigcup_{x \in V(G)} S_x$. Clearly, $S$ is a dominating set of $G \circ H$ since each $S_x$ is a dominating set of $H^x$, for every $x \in V(G)$.

Let $u, z \in S$. To show that $S$ is a 2-movable dominating set of $G \circ H$, we consider the following cases:

**Case 1.** Suppose that there exists $w \in V(G)$ such that $u, z \in S_w$.

Let $a_u \in N_{H^w}(u)$. Now, $(S_w \setminus \{u, z\}) \cup \{a_u, w\}$ is a dominating set of $H^w$. Hence,

$$(S \setminus \{u, z\}) \cup \{a_u, w\} = \left(\bigcup_{x \in V(G) \setminus \{w\}} S_x\right) \cup (S_w \setminus \{u, z\}) \cup \{a_u, w\}$$

is a dominating set of $G \circ H$ since $S_x$ is a dominating set of every $H^x$ for all $x \in V(G) \setminus \{w\}$.

**Case 2.** Suppose that there exist $w_1, w_2 \in V(G)$ such that $u \in S_{w_1}$ and $z \in S_{w_2}$.

Observe that

$$(S \setminus \{u, z\}) \cup \{w_1, w_2\}$$

$$= \left(\bigcup_{x \neq w_1, w_2} S_x\right) \cup ((S_{w_1} \setminus \{u\}) \cup \{w_1\}) \cup ((S_{w_2} \setminus \{z\}) \cup \{w_2\})$$

is a dominating set of $G \circ H$ since $S_x$ is a dominating set of every $H^x$ for all $x \neq w_1, w_2$ and $((S_{w_1} \setminus \{u\}) \cup \{w_1\})$ and $((S_{w_2} \setminus \{z\}) \cup \{w_2\})$ are dominating sets of $H^{w_1}$ and $H^{w_2}$, respectively.

Therefore, $S$ is a 2-movable dominating set of $G \circ H$ and accordingly,

$$\gamma_m^2(G \circ H) \leq |S| = \sum_{x \in V(G)} |S_x| = \sum_{x \in V(G)} \gamma(H) = |V(G)| \gamma(H).$$

Now, suppose that $\gamma_m^2(G \circ H) < |V(G)| \gamma(H)$. Then there exists a minimum 2-movable dominating set $T$ of $G \circ H$ such that $|T| < |V(G)| \gamma(H)$.

Suppose $\gamma(H) = 1$. Then $|T| < |V(G)|$. It follows that $\exists x \in V(G)$ such that $T \cap (V(x + H^x)) = \emptyset$. This contradicts the fact that $T$ is a dominating set of $G \circ H$.

In addition, suppose $\gamma(H) \geq 2$. Since $|T| < |V(G)|\gamma(H)$, $\exists a \in V(G)$ such that $|T \cap V(a + H^a)| < \gamma(H) = \gamma(H^a)$. Let $T_a = T \cap V(a + H^a)$ and $S_a = T \cap V(H^a)$. If $S_a = T_a$, then by Lemma 3.4, $S_a$ is a dominating set of $H^a$. Now, $|T_a| < \gamma(H^a)$. This means that $|S_a| < \gamma(H^a)$. This is a contradiction, since $\gamma(H^a)$ is the minimum cardinality of a dominating set in $H^a$.

On the other hand, suppose $T_a = S_a \cup \{a\}$ and $u \in T_a$. Since $T$ is a 2-movable dominating set, by Lemma 3.5, we have the following:

(i) $T_a \setminus \{a, u\}$ is a dominating set of $H^a$, and

(ii) $((T_a) \setminus \{a, u\}) \cup \{x_a, x_u\}$ is a dominating set of $H^a$ for some $x_a$ and $x_u$ neighbors of $a$ and $u$, respectively.

Now, $|T_a \setminus \{a, u\}| < |T_a| < \gamma(H^a)$ and $|((T_a) \setminus \{a, u\}) \cup \{x_a, x_u\}| = |T_a| < \gamma(H^a)$. This is a contradiction, since $\gamma(H^a)$ is the minimum cardinality of a dominating set in $H^a$.

Thus, $\gamma_m^2(G \circ H) \not< |V(G)|\gamma(H)$. Therefore,

$$\gamma_m^2(G \circ H) = |V(G)|\gamma(H). \qquad \square$$

The next corollary follows from the fact that $K_1 + H \cong K_1 \circ H$.

**Corollary 3.1.** *Let $H$ be a connected graph of order at least* 4. *Then* $\gamma_m^2(K_1 + H) = \gamma(H)$.

**Acknowledgement**

The authors would like to acknowledge the support of the University of Southeastern Philippines in the publication of this research work.

The authors are deeply thankful to the reviewers for their valuable suggestions to improve the quality and presentation of the paper.

**References**

[1]   J. Blair, R. Gera and S. Horton, Movable dominating sensor sets in networks, J. Combin. Math. Combin. Comput. 77 (2011), 103-123.

[2]   O. Ore, Theory of Graphs, American Mathematical Society Colloquium Publications, Vol. XXXVIII, 1962.